# Optimisation in some Banach Algebras related to the Fourier Algebra.

## Edmond E. Granirer

**ABSTRACT.** Let $A_p(G)$ denote the Figa-Talamanca-Herz Banach Algebra of the locally compact group G, thus $A_2(G)$ is the Fourier Algebra of G. If G is commutative then $A_2(G)=L^1(\hat{G})^\wedge$. Let $A_p^r(G)=A_p\cap L^r(G)$ with norm $\|u\|_{A_p^r}=\|u\|_{A_p}+\|u\|_{L^r}$. We investigate a property which insures not only existence of solutions to optimization problems but moreover, facility in testing that an algorithm converges to such solutions, namely the RNP. Theorem(a): If G is weakly amenable then $A_p^r$ is a dual Banach space with the RNP if $1\leq r\leq p'$. This does not hold if G=SL(2,R), p = 2 and r > 2. Theorem(b): If G is weakly amenable and second countable and $A_p^t$ has the RNP for t = s, then it has the RNP for all $1\leq t\leq s$, where $s=\infty$ is allowed. In particular second countable noncompact groups G, for which $A_p(G)$ has the RNP, namely Fell groups, have to satisfy that $A_p^r(G)$ has the RNP for all $1\leq r<\infty$. The results are new, even if G = Z, the additive integers.

## 1. INTRODUCTION

Let G be a locally compact group and let $A_p(G)$ denote the Figa-Talamanca - Herz Banach algebra of G, as defined in [Hz1], thus generated by $L^{p'} * L^{\vee p}(G)$, where $1<p<\infty$, $and\ 1/p+1/p'=1$, see sequel.

Hence $A_2(G)$ is the Fourier algebra of G as defined and studied by Eymard in [Ey1]. If G is abelian then $A_2(G)=L^1(\hat{G})^\wedge$.

Denote $A_p^r(G)=A_p\cap L^r(G)$, for $1\leq r\leq\infty$, equipped with the norm

$$\|u\|_{A_p^r}=\|u\|_{A_p}+\|u\|_{L^r}. \quad \text{If } r=\infty \text{ let } A_p^r(G)=A_p(G).$$

If G is abelian then, $A_2^r(G)=\{u\in L^1(\hat{G})^\wedge; \hat{f}\in L^r(G)\}$, with the norm

$$\|u\|=\|f\|_{L^1(\hat{G})}+\|\hat{f}\|_{L^r(G)}, \quad \text{if } u=\hat{f}.$$

The study of these Banach Algebras started in a beautiful paper of Larsen Liu and Wang [LLW] in the abelian case, and continued in [La1], [La2], [LCh], etc.

Let X be a Banach space. Then

X has the *Krein-Milman Property (**KMP**)*, if *each norm closed convex bounded subset B of X is the norm closed convex hull of it's extreme points, ext(B)*, (namely, $B = \overline{co}\,ext(B)$)

Hence : "O*ptimisation problems in such sets B have solutions*".

The closed unit ball $\mathcal{B}_1$ of $L^1(\mu)$, for a nonatomic measure $\mu$, has no extreme points, hence it *does not* have the **KMP.** And yet, if $\mu$ is atomic (for example in case of $\ell^1$) then $\mathcal{B}_1$ has many extreme points. In fact this space does even have the **KMP** and is in addition a dual space. Known results (see [DU1]) then imply that this space has a stronger property denoted **RNP**, namely:

X has the *Radon-Nikodym Property (**RNP**)* if *each B as above is the norm closed convex hull of it's strongly exposed points (strexp(B))*, see sequel or [DU1], p. 190 and p. 218.

Points in *strexp(B)* are points of *ext(B)* that have beautiful smoothness properties. In particular they are weak to norm continuity points of B and are peak points of B.

Hence: "*Optimisation problems in such sets B have solutions, but moreover, it is easy to test if an algorithm converges to a solution*".

Quoting Jerry Uhl: "*A Banach space has the **RNP** if it's unit ball wants to be weakly compact, but just cannot make it*".

**Definition:** *Let B be a bounded subset of the Banach space X and $b \in B$. b is a strongly exposed point of B (and strexp(B) denotes the set of all such), if*
$\exists b^* \in X^*$ *such that*

$$\operatorname{Re} b^*(x) < \operatorname{Re} b^*(b), \forall x \in B \text{ and } x \neq b, \text{ and}$$

$\operatorname{Re} b^*(x_n) \to \operatorname{Re} b^*(b) \text{ for } x_n \in B$ *implies* $\|x_n - b\| \to 0$. *(see [DU1] p.138)*

*Hence in order to test an algorithm for* $b \in str\exp B$, *it is enough to test it on one element of* $X^*$.

Any X which is norm isomorphic to $\ell^1$ has the **RNP**. If X a dual Banach space, and B is w* compact convex, then the functional b* can be chosen in the predual of X, see [DU1] p.213.

It follows from above that if G is *abelian* then $A_2(G)$ *has the **RNP** if G is compact and does not have the **RNP** it if G is not compact.*

And yet, for any abelian G and any compact subset K,
$A_K^2(G) = \{u \in A_2(G): spt\,u \subset K\}$, *does have the **RNP***, where spt denotes support.

In fact we have proved in [Gr1] that for *any* G and any compact subset K and any $1<p<\infty$, $A_K^p = \{u \in A_p(G); spt\, u \subset K\}$ has the **RNP**. Tools for abelian G are not available in this case.

It has been proved by W. Braun, in an unpublished preprint [Br], that if *G is amenable, then $A_p^1(G)$ is a dual Banach space with the RNP*. The result in [Br] uses the method in [Gr1] and the involved machinery of [BrF], which is avoided below.

Denote as in [Hz1]

$$A_p(G) = \{u = \sum u_n * v_n^{\vee}; u_n \in L^{p'}, v_n \in L^p, \sum \|u_n\|_{L^{p'}} \|v_n\|_{L^p} < \infty\},$$

where the norm of $u \in A_p$ is the infimum of the last sum over all the representations of u as as above.

We will omit at times G and write $L^p, A_p$, etc. instead of $L^p(G), A_p(G)$, etc.

Denote by $PM_p(G) = A_p(G)^*$, and by $PF_p(G)$, the norm closure in $PM_p(G)$ of $L^1(G)$, ( as a space of convolutors on $L_p(G)$).

Let $W_p(G) = PF_p(G)^*$. Then $W_p(G)$ is a Banach algebra of bounded continuous functions on G containing the ideal $A_p(G)$, studied by Cowling in [Co1]. Let $W_p^r(G) = W_p \cap L^r(G)$

If G is abelian and p=2 then $W_2(G) = M(\hat{G})^{\wedge}$, where *M(G)* is the space of bounded Borel measures on G. Let $C_0(G)[C_c(G)]$ denote the continuous functions which tend to 0 at $\infty$, [with compact support], with norm $\|u\|_{\infty} = \sup\{|u(x)|; x \in G\}$.

The *group G is weakly amenable* if $A_2(G)$ has an approximate identity $\{v_{\alpha}\}$ bounded in the norm of $B_2(G)$, the space of Herz-Schur multipliers, see [Hz1]-[Hz2], [Ey2] (or [DCH], [Gr5]). As known the free group on n > 1 generators is weakly amenable but non amenable. For much more see [DCH].

Our first result hereby is the

**Theorem1:** *Let p=2 or G be weakly amenable and $1<p<\infty$.*
*Then*

$$(*)\, W_p \cap L^r(G) = A_p \cap L^r(G), \forall 1 \leq r \leq p'.$$

*Hence ( by [Gr5] Thm. 2.2) $A_p^r(G)$ is a dual Banach space for such r.*

*If G is in addition unimodular then the above holds for $1 \leq r \leq \max(p, p')$.*

The interval [1, p'] cannot be improved even if p=2 and G = Z, the additive integers, see sequel.

Use of Theorem 1 is made in proving the main result of this section, namely:

**Theorem 2:** *Let G be a weakly amenable locally compact group. Then*
$\forall 1 \le r \le p'$, $A_p^r(G) = W_p^r(G)$ *and* $A_p^r(G)$ *is a dual Banach Algebra with the* **RNP**.
*If G is in addition unimodular this is the case* $\forall 1 \le r \le \max(p, p')$.

**Remark:** If G=SL(2,R), p = p'= 2 and r > 2, then $A_2^r(G) = A_2(G)$ (see [KuS], [Co2]) and $A_2^r(G)$ does not have the **RNP** and is not a dual space, see [Gr4] p. 4382. Note that G is a weakly amenable, but non amenable group.
The unimodular case, for both the above results, was proved in our paper [Gr5].

**Remarks:** **(1)** A group G with completely reducible regular representation is called in [T] an *[AR] group*. G is such iff $A_2(G)$ has the **RNP**, as proved by Keith Taylor [T]. (see[T] for much more). A noncompact [AR] group is called a Fell group, see [B] section IV.
Larry Baggett and Keith Taylor construct in [BT] p.596 (iii) an example of a connected nonunimodular Lie group G = G₃ such that $A_2(G) \ne W_2 \cap C_0(G)$ and such that G is a Fell group. Our next result will imply that for this group G, $A_2^r(G)$ *has the* **RNP**, *for all r*.

If [BT] could be improved to show that for some finite s > 2, $A_2^s(G) \ne W_2^s(G)$, it would follow that $A_2^s(G)$ having the **RNP** *does not imply that* $A_2^s(G) = W_2^s(G)$.
(Necessarily s > 2, if G = G₄ in [BT] p.597, (which is amenable) since then $A_2^r(G) = W_2^r(G), \forall 1 \le r \le 2$, by **Theorem 2**).

**(2)** Assume that for arbitrary G, $A_2^s(G)$ having the **RNP** for s > 2 would imply the equality $W_2^s(G) = A_2^s(G)$. It *would then follow for G = Z, that* $A_2^s(Z)$ *does not have the* **RNP** *for all s > 2. This is implied by the fact that* $A_2^s(Z) \ne W_2^s(Z), \forall s > 2$ *by the Hewitt-Zuckerman [HZ] result, noted in [LiR] and used in [Gr4] p. 4379. Hence there would be no need to take G = SL(2,R) in the remark above, and Z would suffice.*

In the next section we are interested in the following question:

*Given p and the group G, determine the set O(p,G) of those r for which* $A_p^r(G)$ *has the RNP, ( O for optimization!)*
We will show, using results on semi embeddings due to H.P.Rosenthal ,(see [R], [LPP]) the following

**Theorem 3 :** *Let $G$ be second countable and $A_p(G)$ have a multiplier bounded approximate identity: If $A_p^t(G)$ has the RNP for $t = s$ then $A_p^t(G)$ it has the RNP for all $1 \leq t \leq s$, where $s = \infty$ is allowed.*

The above results show that:
 **(1)** $[1,p'] \subset O(p,G)$, *if $G$ is weakly amenable,*
 and $[1,\max(p,p')] \subset O(p,G)$ *if $G$ is in addition unimodular. Note that weak amenability depends only on 2 , yet the result holds for all p .*
 **(2)** $[1,2] = O(2,G)$, *if $G = SL(2,R)$. This shows that **(1)** is the best one can do.*
 **(3)** $[1,\infty] = O(2,G)$, *if $G$ is a Fell group.*

2. MAIN RESULTS

(I) **NO UNIMODULARITY**. We improve hereby results in [Gr4], [Gr5], by removing the unimodularity of the group G.

**Theorem1 :** *Let G be a locally compact group. If $p=2$,*

*or if G is weakly amenable and $1 < p < \infty$ , then*

$$(*) \quad W_p \cap L^r(G) = A_p \cap L^r(G), \forall 1 \leq r \leq p'.$$

*If G is in addition unimodular then (\*) holds for $\forall 1 \leq r \leq \max(p, p')$.*

*Hence , for the above values of r, $A_p^r(G)$ is a dual Banach space.*

 **Remark** *The interval $[1,p']$ is the best one can do even for $G=Z$ and $p = 2$ as proved by Hewitt and Zuckerman in [HZ], and as noted in [LiR]. (see also [Gr4] p.4379.)*
 *If $G=SL(2,R)$ or $SL(2,C)$ and $p=2$ then (\*) does not hold for any $r > 2$, and in addition $A_p^r$ is not a dual space for $r > 2$.*

**Proof:** By weak amenability, for all $1 < p < \infty$ , the $W_p$ norm restricted to $A_p$ is equivalent to the $A_p$ norm , ([Gr5] Corollary 3.7.) . If $p = 2$ then Kaplansky's density theorem will yield the same result.

Now with the notations of [Gr 4] Thm. 2.1. p. 4379 , if $e_\alpha \in C_c(G)$ is an approximate identity for $L_1(G)$ , such that each $e_\alpha$ is the square of a special operator , a la Fendler [Fe] p.129 , we have , loc. cit. $\|e_\alpha * w - w\| \mapsto 0 \ \forall w \in W_{p'}$ afortiori $\forall w \in W_{p'} \cap L^{p'\vee}$.

But, since $e_\alpha \in C_c(G)$, we have for such $w$, that

$e_\alpha * w \in L^p * L^{p'\vee} \subset A_{p'}$, thus $e_\alpha * w$ is a Cauchy sequence in $A_{p'}$.

Hence $w \in A_{p'}$. It follows that $W_{p'} \cap L^{p'\vee} = A_{p'} \cap L^{p'\vee}$.

However by [Co] p.91, $W_p = W_{p'}{}^\vee$, $A_p = A_{p'}{}^\vee$. Hence
$$W_p \cap L^{p'} = A_p \cap L^{p'}, \forall 1 < p < \infty.$$

But $W_p$ contains only bounded functions, hence

$\forall r \leq p'$, $W_p \cap L^r = W_p \cap L^{p'} \cap L^r = A_p \cap L^{p'} \cap L^r = A_p \cap L^r$. Thus

(i) $W_p \cap L^r = A_p \cap L^r, \forall r \leq p'$.

If $G$ is unimodular then, since $(W_p \cap L^r)^\vee = (A_p \cap L^r)^\vee$,

it follows that

(ii) $W_{p'} \cap L^r = A_{p'} \cap L^r, \forall r \leq p'$, which holds for all $1 < p' < \infty$.

Replace now p' by p in (ii) and then

(iii) $W_p \cap L^r = A_p \cap L^r, \forall r \leq p$.

Thus (i) and (iii) imply the unimodular case.

By Theorem 2.2 of [Gr5] $W_p(G) \cap L^r(G)$ is a dual Banach space for all $1 < p < \infty$ and $1 \leq r \leq \infty$, and all locally compact groups $G$. □

**Lemma1:** *Let G be a locally compact group. Assume that $A_p(G)$ has an approximate identity $u_\alpha$ such that $\sup \|u_\alpha\|_\infty \leq B < \infty$. Then*

*(a) $A_p \cap C_c$ is norm dense in $A_p^r$ and*

*(b) If G is second countable then $A_p^r$ is norm separable.*

**Proof:** (a) Let $e_\alpha \in A_p \cap C_c$ satisfy $\|e_\alpha - u_\alpha\|_{A_p} \to 0$ and $\|e_\alpha - u_\alpha\|_{A_p} \leq 1, \forall \alpha$.

Then $\|e_\alpha\|_\infty \leq \|e_\alpha - u_\alpha\|_\infty + \|u_\alpha\|_\infty \leq 1+B$. Hence

$$\|e_\alpha v - v\|_{A_p} \leq \|(e_\alpha - u_\alpha)v\|_{A_p} + \|u_\alpha v - v\|_{A_p} \to 0, \forall v \in A_p.$$ But if $w \in A_p^r$ and

$K \subset G$ is compact such that $\int_{G \sim K} |w|^r dx < \epsilon$ then

$\int_{G \sim K} |(e_\alpha - 1)w|^r \leq \int_{G \sim K} (2+B)|w|^r \leq (2+B)\epsilon$. But $\int_K |(e_\alpha - 1)w|^r \to \infty$. It thus follows

that $\|e_\alpha w - w\|_{A_p^r} \to 0$. But $e_\alpha w \in A_p \cap C_c$.

(b) $A_p(G)$ is norm separable, hence so is $A_p[K] = \{u \in A_p(G); spt\, u \subset K\}$, where

$K \subset G$. Let $A_p^r[K] = \{u \in A_p^r(G); spt\, u \subset K\}$. If K is compact then the identity

$I: A_p^r[K] \to A_p[K]$ is 1-1, onto and continuous, hence it is bi continuous. Hence

$A_p^r[K]$ is separable. Let now $K_n \subset int\, K_{n+1} \subset G$, be compact (int denotes interior),

such that $\cup K_n = G$. It is hence enough to show that $\cup A_p^r[K_n]$ is norm dense in $A_p^r(G)$

By (a) we know that $A_p \cap C_c$ is norm dense in $A_p^r(G)$. But if $v \in A_p^r(G)$ has compact

support S then $S \subset K_j$ for some j, hence $v \in A_p^r[K_j]$. Thus $\cup A_p^r[K_n]$ is norm dense

in $A_p^r(G)$. □

*Remark:* We do not know if, $A_p \cap C_c(G)$ is norm dense in $A_p^r(G)$ even for

$G = SL(2,R) \triangleleft R^2$, if $p = 2$ and all r, (see [Do]).

*Corollary1:* Let G be a second countable locally compact group. If G is weakly
amenable then $\forall 1 \leq r \leq p'$, $A_p^r(G)$ is a separable dual Banach algebra and thus has
the RNP.
 If G is in addition unimodular, this is the case $\forall 1 \leq r \leq \max(p, p')$.

*Remark: Weak amenability, namely the existence in $A_2$ of an approximate identity norm
bounded in $B_2$ depends only on p = 2, yet the result holds for all p. Since by Furuta's
Thm.2.4 in [Fu], $B_2 \subset B_p$ contractively, see also [Gr5] p.23.*

*Proof:* By Thm. 2.2 on of [Gr3] and the above Proposition $A_p^r(G)$ is a dual Banach
space $\forall 1 \leq r \leq p'$. But since G is weakly amenable $A_p(G)$, has a multiplier bounded approximate identity $\forall 1 < p < \infty$, by [Gr2013] p.24. It thus follows by the Lemma above, that $A_p^r(G)$ is norm separable. But separable dual Banach spaces have the RNP by [DU] p.218. □

The main result of this section is the

**Theorem2:** *Let G be a weakly amenable locally compact group. Then*
$\forall 1 \leq r \leq p'$, $A_p^r(G) = W_p^r(G)$ *and* $A_p^r(G)$ *is a dual Banach algebra with the RNP.*
  *If G is unimodular, this is the case* $\forall 1 \leq r \leq \max(p,p')$.

**Proof:** Based on the above Corollary follow the proof of Theorem 3.1 on p.p.22-24 of [Gr2013] and [Gr2011P] p.4381. □

**Remark:** If G=SL(2,R), p=p'=2 and r > 2 then $A_2^r(G) = A_2(G)$ and $A_2^r(G)$ does not have the RNP and is not a dual space, see [Gr4] p.4382.

**(II) INTERVALS WITH THE RNP.** We will show that if G is second countable and $A_p(G)$ has a multiplier bounded approximate identity then $A_p^t(G)$ having the RNP for t = s implies that it has it for all $1 \leq t \leq s$, where $s = \infty$ is allowed.

**Definition:** *Let X,Y be Banach spaces and* $T: X \to Y$ *be a bounded linear operator. T is a semi-embedding if it is **one to one** and **it maps the closed unit ball in X into a closed set in Y**. If such T exists we say that X semi embeds in Y.*

**Theorem [H.P.Rosenthal]:** *A separable Banach space has the RNP if it semi-embeds in a Banach space with the RNP.* See [DU2] p.160 or [Ro], [LPP].

We will make use of the above Theorem of Rosenthal, to prove the main Theorem3.

We will need the following

**Lemma2:** *If r < s then the identity* $I : A_p^r(G) \to A_p^s(G)$ *is a semi-embeding, for any* $s \leq \infty$ *(If* $s = \infty$, $A_p^\infty(G) = A_p(G)$).

**Proof:** Denote by $B_t$ the closed unit ball of $A_p^t(G)$. Let $v_n \in B_r$ satisfy that $\|v_n - w\|_{A_p^s} = \|v_n - w\|_{A_p} + \|v_n - w\|_{L^s} \to 0$, for some $w \in A_p^s(G)$.

(If $s = \infty$ only $\|v_n - w\|_{A_p}$ appears). Clearly $|v_n(x)| \to |w(x)|, \forall x \in G$. And by Fatou's Lemma we have $\int |w|^r dx \leq \liminf \int |v_n|^r dx \leq 1$ Thus $w \in A_p^r$. But
$1 \geq \limsup(\|v_n\|_{A_p} + \|v_n\|_{L^r}) \geq \lim \|v_n\|_{A_p} + \liminf \|v_n\|_{L^r} \geq \|w\|_{A_p} + \|w\|_{L^r}$. Thus $w \in B_r$. □

**Theorem 3:** Assume that G is second countable and $A_p(G)$ has a multiplier bounded approximate identity.
If for some $t \leq \infty$, $A_p^t(G)$ has the RNP, then so does $A_p^r(G), \forall 1 \leq r \leq t$.
In particular, if $A_p(G)$ has the RNP then $A_p^r(G)$ has the RNP for all $1 \leq r < \infty$.

**Proof:** Apply Rosenthal's Theorem and the above Lemma 2, and note that, by Lemma 1 $A_p^r(G)$ is necessarily norm separable, due to the existence of the multiplier bounded approximate identity (thus bounded in the uniform norm). □

**Remark:** *Note that by [Fu] Thm.2.4 and p.581, weak amenability implies existence of a multiplier bounded approximate identity.*

*Note that* the Fell group in [B] p.142 is unimodular and CCR.

Mauceri and Picardello have constructed in [MP], amenable and nonamenable, totally disconnected unimodular Fell groups, generalizing the original Fell group. Many of these are p-adic matrix groups.

Baggett and Taylor present in [BT] examples of Fell groups which are connected Lie groups and which are (i) solvable, (ii) amenable nonsolvable, (iii) nonamenable, (iv) non-TypeI. All of which are not unimodular.

For all the above ones, $A_2^r(G)$, $\forall 1 \leq r \leq \infty$ has the RNP, by Rosenthal's theorem and the above lemma.

**Question:** *If G is noncompact abelian then $A_2(G)$ **does not** have the RNP, since its closed unit ball has no extreme points see [DU] p.219. Yet, $A_2^r(G), \forall 1 \leq r \leq 2$ **does have the RNP**, by Theorem 1. For such G nothing is known if r > 2.*

**REFERENCES**

[B]  L. Baggett, A separable group having discrete dual is compact. J. Functional Anal. 10 (1972), 131-148.
[BT]  Larry Baggett and Keith Taylor, Groups with Completely Reducible Regular Representation. Proc. Amer. Math. Soc. 72 (1978), 593-600.
[Br]  W. Braun: Einige Bemerkungen Zu S_0(G) und A^p(G)intL^1(G). Preprint.
[BrF]  W. Braun and Hans G. Feichtinger. Banach Spaces of Distributions Having Two Module Structures. J.Funct. Analysis. 51 (1983) 174-212.
[Co1]  Michael Cowling, An Application of Littlewood-Paley Theory in Harmonic Analysis. Math. Ann. 241 1979), 83-96.

[Co2]   Michael Cowling,  The Kunze-Stein phenomenon. Ann.Of Math.  106 (1978), 209-234.
[DCH] J. deCanniere and U. Haagerup: multipliers of the Fourier algebra of some simple Lie groups and their discrete subgroups. Amer. J. Math. 107 (1985), 455–500.
[DU1]  J.Diestel and J.J.Uhl,Jr: Vector Measures. Math. Surveys. Amer. Math. Soc. 1977.
[DU2]  J. Diestel and J.J. Uhl,Jr : Progress in Vector measures – 1977-83. Measure Theory and Applications. LNM. 1033. Springer. 1983
[Do]  B. Dorofaeff : The Fourier Algebra of  $SL(2,R)/\backslash R_n$ , n>1 , has no multiplier bounded approximate unit. Math. Ann. 297 (1993, 707-724.
[Ey1] P. Eymard: L'algebre de Fourier d'un groupe localement compacte. Bull. Soc. Math. France. 92 (1964), 181–236.
[Ey2] P. Eymard: Algebre $A_p$ et convoluteurs de $L_p$ . Lecture Notes in Math. No. 180 (Springer 1971), 364–381.

[Fe]  Gero Fendler : An Lp version of a theorem of D. A. Raikov. Ann. Inst. Fourier, Grenoble. 35 (1985), 125-135.
[Fu]  Koji Furuta :  Algebras  $A_p$  and  $B_p$  and amenability of locally compact groups. Hokkaido Math. J. 20 (1991) 579-591.
[Gr1]  Edmond E. Granirer : An Application of the Radon Nikodym Property in Harmonic   Analysis. Bull. U.M.I. (5) 18-B (1981), 663-671.
 [Gr2]:     -------   :Amenability and semisimplicity for second duals of quotients of the Fourier Algebra A(G).  J. Austral. Math. Soc. (Series A) 63 (1997), 289-296.
[Gr3]     -------      :  The Figa-Talamanca-Herz-Lebesgue Banach Algebras
$A_p^r(G) = A_p \cap L^r(G)$ . Math. Proc. Camb. Phil. Soc. 140 (2006), 401-416.
[Gr4]     -------      : The Radon-Nikodym Property for some Banach Algebras related to the Fourier Algebra. Proc. Amer. Math. Soc. 139 (2011), 4377-4384.
 [Gr5]     --------      : Weakly Amenable Groups and the RNP for some Banach Algebras related to the Fourier Algebras.  Coll. Math. 130 (2013), 19-26.
[HZ]  Edwin Hewitt and Herbert Zuckerman : Singular measures with absolutely continuous convolution squares. Proc. Camb. Phil. Soc. 62 (1966), 399-420.
[Hz1]  C. Herz: Harmonic Synthesis for Subgroups. Ann. Inst. Fourier, Grenoble. 23 (1973) 91-123.
[Hz2] C. Herz: The theory of p spaces with an application to convolution operators. Trans. Amer. Math. Soc. 154 (1971) 69–82.
[LLW]  Ron Larsen, Ten-sun Liu, Ju-kwei Wang: On functions withFourierTransforms in $L_p$ . Mich. Math. J. 11 (1964), 369-378.
[La1]  Hang-Chin Lai: On some properties of  A^p(G) algebras. Proc.Japan Acad. 45 (1969), 572-576.
[La2]  Hang-Chin Lai: A remark on A^p(G) algebras. Proc. Japan Acad. 46 (1970) 58-  63.
[LPP]  Lotz.H.P, Peck.N.T., Porta. H. Semi-embeddings of Banach Spaces. Proc.Edinburgh Math. Soc. 22 (1979) , 233-240
 [LiR]   Teng-sun Liu and Arnoud van Rooij : Sums and Intersections of Normed Linear Spaces. Math. Nachrichten. 42 (1969), 29-42.

[MP] G. Mauceri and M.A. Picardello : Noncompact unimodular groups with purely atomic Plancherel measures . Proc. Amer. Math. Soc.78 (1980) , 77-84.

[Ri] N.W.Rickert : Convolutions of Lp functions. Proc. Amer. Math. Soc. 18 (1967), 762-763. MR0216301 (35:7136)

[R] Rosenthal H.P. Convolution by a Biased Coin. The Altgelt Book 1975/76. University of Illinois Functional Analysis Seminar.

[Sa] Sadahiro Saeki : The Lp conjecture and Young's Inequality. Ill. J. Math. 34 (1990), 614-627.

[T] Keith Taylor, Geometry of the Fourier Algebras and Locally Compact Groups with Atomic Unitary Representations. Math. Ann. 262 (1983), 183-190.

Dept. of Math. Univ. of B.C. Vancouver B.C. V6T IZ4, Canada.
E-mail address: granirer@math.ubc.ca